\documentclass[11pt]{amsart}
  \makeatletter
\def\@normalsize{\@setsize\normalsize{10pt}\xpt\@xpt
\abovedlayskip 10pt plus2pt minus5pt\belowdisplayskip \abovedisplayskip \abovedisplayshortskip \z@
plus3pt\belowdisplayshortskip 6pt plus3pt minus3pt\let\@listi\@listI}
\def\subsize{\@setsize\subsize{12pt}\xipt\@xipt}
\def\section{\@startsection {section}{1}{\z@}{1.0ex plus 1ex minus
2ex}{.2ex plus .2ex}{\large\bf}}
\def\subsection{\@startsection {subsection}{2}{\z@}{.2ex plus 1ex}
{.2ex plus .2ex}{\subsize\bf}} \makeatother

\newtheorem{theorem}{Theorem}

\theoremstyle{remark}
\newtheorem*{remark}{Remark}

\begin{document}

\baselineskip=12pt

\address{  Helmut Kaeutner Str. 25, 81739 Munich, Germany, tel. 0896256676 }



\email{alexey.ostrovskiy@unibw.de}



\title{\bf {$\sigma$-homogeneity of Borel sets
 }}
\author {Alexey Ostrovsky}
 \maketitle

\begin{abstract}

We give an affirmative answer to the following question: 

Is any Borel subset of a Cantor set $\textbf{ C}$ a sum of a
countable number of pairwise disjoint   
$h$-homogeneous subspaces that
  are closed in $X$?

It follows that  every  Borel  set   $X \subset  \textbf{ R}^n$
can be partitioned into countably many $h$-homogeneous  subspaces
that are $G_{\delta}$-sets in $X$.

 \end{abstract}

\noindent

\vspace{0.15in}

We will denote by  $\textbf{R}$ ,  $\textbf{P}$,  $\textbf{Q}$,
and $\textbf{C}$   the spaces of real, irrational, rational
numbers, and a Cantor set, respectively.

\smallskip

 Recall that a zero-dimensional topological space $X$ is
 \emph {$h$-homogeneous} if  $U$ is homeomorphic to $X$ for each
 nonempty clopen subset $U \subset X$.  More about topological properties  of $h$-homogeneous spaces  see, for example,  in
 \cite{Ma}, \cite{Me},  \cite{Mo}, \cite{Te}.

 \smallskip

  \vspace{0.10in}

    We call a  zero-dimensional metric space $X$  \emph{$\sigma$-homogeneous}
    if it  is a countable union  of $h$-homogeneous subspaces $X_i$ that are closed in $X$.
    It is easily seen that every set  $X_i  \setminus  \bigcup_{j <i}X_j$ is an open
    subspace in $X_i$  and   can be partitioned into countably many pairwise disjoint  subsets that are clopen in $X_i$.

  \smallskip

   Hence, a space $X \subset \textbf{ C}$ is
  $\sigma$-homogeneous iff it can be partitioned into a countably many  pairwise disjoint $h$-homogeneous subspaces that
  are closed in   $X$.

According to  the Cantor--Bendixson theorem,  every closed subset
$F \subset \textbf{C}$  is $\sigma$-homogeneous.

\smallskip

The question of whether this  assertion holds  for all  Borel
subsets of $\textbf{C}$ was posed  in  \cite [p.228] {BAO}.

         \vspace{0.10in}

The following  theorem  gives an affirmative answer to the above question.

  \vspace{0.15in}

\begin{theorem}

{Every Borel set  $ X \subset \textbf{C}$ is a $\sigma$-homogeneous space. }

\end{theorem}

  \vspace{0.10in}

\textbf{Proof}. We  proved the following simple proposition in
   \cite [Theorem 7] {OA}:

\smallskip

Every  $\Pi_2^0$-set (and every  $\Sigma_2^0$-set)
     $X \subset \textbf{C}$ is representable  as a union of
     countably many disjoint closed  copies of following spaces:

\smallskip

     (a)  a singleton set;

\smallskip

      (b) a Cantor set \textbf{C};

      \smallskip

     (c) irrational numbers \textbf{P}.

\smallskip

From the topological characterization of  \textbf{C} and
\textbf{P} it follows that they (and  obviously a singleton set)
are $h$-homogeneous.

\smallskip

 Recall (for more detail  we refer the reader to  \cite{vE}) that $ A \leq_{w} B$ 
if for some continuous $f :  \textbf{C} \to  \textbf{C}$ we have $A = f^{-1}(B) $.

\smallskip

The Borel Wadge  class of a Borel set $A$ is $[A]= \{B \subset  \textbf{C}: B   \leq_{w} A \}.$

\smallskip

The Wadge ordering $<$ on dual pairs $\{\Gamma, \check\Gamma
\}$ (where $ \check\Gamma = \{\textbf{C} \setminus A: A \in
\Gamma  \}$) of Wadge classes that   well-orders the pairs of
Borel Wadge classes is defined by

\smallskip

$\{\Gamma_0, \check\Gamma_0 \}$ $<$ $\{\Gamma_1, \check\Gamma_1
\}$ if and only if $\Gamma_0 \subset \Gamma_1$ and $\Gamma_0 \not
= \Gamma_1.$

 $\Gamma$ is self-dual if $\Gamma = \check \Gamma$.


\smallskip

Also,  for the  classes $\Gamma_{u_0}$ and $\check\Gamma_{u_0}
$, where $\Gamma_{u_0}$ is the class of $F_{\sigma}$-sets and
$\check \Gamma_{u_0}$ is the class of $G_{\delta}$-sets,  
Theorem 1 was proved in  \cite{OA}.

\smallskip

We make an induction hypothesis that the theorem is valid for all
$\Gamma_{\alpha}$ and $ \check \Gamma_{\alpha}$ for all ${\alpha}
< {\beta}$.

\smallskip

Below, we consider two cases \textbf{1} and \textbf{2}.

 \smallskip

 \textbf{1. }Suppose  $\Gamma_{\beta}$ is not a self-dual class and
$X \in \Gamma_{\beta} \setminus \check \Gamma_{\beta}$ .

\smallskip

\textbf{1.1. } If $X$ contains a clopen set $U_1$ of some
class $\Gamma_{\alpha_1}$, ${\alpha_1} < {\beta}$, then $U_1$
falls under the induction hypothesis, and we then consider the set
$X_1 = X \setminus U_1$. Obviously, $X_1 $ is closed in $X$.

\smallskip

If $X_1$  contains a clopen (in $X_1$) subset $U_2$  of some
class $\Gamma_{\alpha_2}$, ${\alpha_2} < {\beta}$, then it falls
under the induction hypothesis, and we then consider the closed
set $X_2 = X_1 \setminus U_2$.

Continuing this process as above, we  get a chain of closed sets
$$X \supset X_1 \supset ...\supset X_{\gamma} \supset...$$
($X_{\gamma} = \bigcap_{\beta < \gamma} X_{\beta}$ for the limit
$\gamma$) that, as we know, stabilizes
  at some countable $\gamma_0  < \omega_1$; i.e., $ X_{\gamma_0} =  X_{\gamma_{0+1}}= ...$.

\smallskip

It is clear that  $X_{\gamma_0}$ is a closed set.

\smallskip

Obviously, $X$ is a countable union of pairwise
   disjoint closed  sets $U_{\alpha}, (\alpha < \gamma_0)$ and $X_{\gamma_0}$.

\smallskip

 If $X_{\gamma_0} \in \Gamma_{\alpha}$ with $\alpha <
\beta$, then 
the theorem is proved since the sets $U_{\alpha}$ and
$X_{\gamma_0}$ fall
under the induction hypothesis. 
Hence we can suppose that $X_{\gamma_0}$ is  nonempty and
everywhere $ \Gamma_{\beta} \setminus \check
\Gamma_{\beta}.$

\smallskip

\textbf{1.2.} If  $X_{\gamma_0}$ is everywhere of the second
category,  then we get the theorem since by theorems
Keldysh, Harrington and Steel    \cite{Ke},  \cite{OA},
\cite{St}  all the spaces  everywhere of the second category and
everywhere   $\Gamma_{\beta} \setminus \check \Gamma_{\beta}$
(for non-$F_{\sigma}$ or non-$G_{\delta}$ classes) are
homeomorphic.

\smallskip

\textbf{1.3.} Let $X_{\gamma_0}$ be not everywhere of the second
category and, hence, contains a clopen (in $X_{\gamma_0}$)
subset $T_1$  of the first category.

\smallskip

If  $Y_1 = Y \setminus T_1$ contains a clopen set $U_1$ of
some class $\Gamma_{\alpha_1}$, ${\alpha_1} < {\beta}$, we can
repeat  the process of \textbf{1.1}, etc.

It is clear that we obtain
by this way  a subspace $T$ that is (everywhere) of the first category
and everywhere $\Gamma_{\beta} \setminus \check
\Gamma_{\beta}$, which is  $h$-homogeneous  by theorems
Keldysh, Harrington and Steel.

\smallskip

\textbf{2.0.} $X \in \Gamma \cap \check  \Gamma$. Then     \cite[Lemma 4.4.1] {vE}
  there is a nonempty clopen subset $D \subset X$
such that
  $  [X \cap D] <_w [X]$ and $X$ is decomposed into sets of lower
  Wadge rank.

We can repeat the process of \textbf{1.1}.

 \qed

\vspace{0.1in}

Since    $\textbf{R} = \textbf{P} \cup \textbf{Q}$ and   \textbf{P}
is homeomorphic to some  $G_{\delta}$ subset in \textbf{C}  and \textbf{R},
we obtain  the following corollary:

\vspace{0.1in}

\emph{ Every   Borel set $X \subset  \textbf{R}^ n$ can be
partitioned into countably many $h$-homogeneous
$G_{\delta}$-subspaces. }

\smallskip

\vspace{0.15in}

Questions  on the number of topological types  of  homogeneous
Borel sets have been posed by Aleksandrov and Urysohn, Lusin, 
Keldysh \cite{AU},  \cite{L}, \cite{Ke}.

\smallskip

By Keldysh's theorem, every Borel set in $\textbf{C}$ is a
countable sum of \emph{canonical elements} that are homeomorphic
to $\textbf{P}$, $\textbf{C},$ a singleton set or $h$-homogeneous
$\Pi_{ \alpha }^0$-sets  (which are not $\Sigma_{ \alpha
}^0$-sets,  $\alpha >1$)  of the first category in themselves.
 \cite{OSTR}
 \cite{Ke}.

  \smallskip

Since  $\textbf{P},$ $\textbf{C}$, and a singleton set are spaces
of the second category in themselves, it would be  reasonable to
find an analogue of Keldysh's theorem  for $h$-homogeneous
subspaces of the second category.  Using the following simple observation (see also \cite{OA})
we show below  that the assertion of  Keldysh's theorem holds for  $h$-homogeneous
subspaces of the second category in themselves.

\begin{remark}

If   $X$ is of the first  category in itself  everywhere   $ \Gamma_{\beta}
\setminus \check \Gamma_{\beta}$,  where  $\Gamma_{\beta}  =\Pi_{
\alpha }^0$ $(\alpha >1$),  then   $ X$ is  homeomorphic  to  the product $Y \times
\textbf{Q}$, where $Y$ is  a space   everywhere    $ \Gamma_{\beta}
\setminus \check \Gamma_{\beta}$  of the second  category in itself.

\end{remark}

 Indeed, denote (all embeddings in \textbf{C} are dense):

 \vspace{0.15in}

 $ \textbf{C}_1 = cl_{\textbf{C}}X$ (it is clear that  $\textbf{C}_1$ is homeomorphic to
$ \textbf{C});$

$Y =\textbf{C} \setminus ((\textbf{C}_1 \setminus X)) \times \textbf{Q})$.

\vspace{0.05in}

Obviously, $Y$ is everywhere  $ \Gamma_{\beta} \setminus \check \Gamma_{\beta}$
 of the second category in itself and $Y \times \textbf{Q}$ is 
everywhere $ \Gamma_{\beta} \setminus \check \Gamma_{\beta}$
of the first category  in itself.

  \vspace{0.15in}
  
 Finally,    $Y \times
\textbf{Q}$ is  homeomorphic  to  $ X$\footnote{Note that     S. Medvedev  proved  that  every $h$-homogeneous space $Y \subset \textbf{C}$
        of the first category in itself is homeomorphic to  $Y \times \textbf{Q}.$ For  more details we refer the reader to  \cite{Me}.        
}. Hence,  every canonical
element of Keldysh $X$  is a sum of a countable number of pairwise
disjoint   $h$-homogeneous subspaces  of the second category
(that
  are closed in $X$).

\end{document}